\UseRawInputEncoding
\documentclass[a4paper,11pt]{article}
\usepackage{graphicx}
\usepackage{amsmath}
\usepackage{amssymb}
\usepackage{amstext}
\usepackage{color}
\usepackage{colortbl}
\usepackage{xcolor}

\newcommand{\proof}{\vskip 5truemm \noindent{\textsc Proof.}~}

\newtheorem{rem}{Remark}[section]
\newtheorem{thm}{Theorem}[section]
\newtheorem{coro}{Corollary}[section]
\newtheorem{lem}{Lemma}[section]

\numberwithin{equation}{section}
\begin{document}

\title{\bf BV right-continuous solutions of second-order differential inclusions governed by maximal monotone operators}

\author{Dalila Azzam-Laouir
 \footnote{LAOTI, FSEI, Universit\'e Mohamed Seddik Benyahia de Jijel, Alg\'erie. E-mail: laouir.dalila@gmail.com,
 dalilalaouir@univ-jijel.dz}}

\maketitle

\begin{abstract} This paper concerns  existence
of right-continuous with bounded variation solutions of a perturbed
second-order differential inclusion governed by time and
state-dependent maximal monotone operators.
\end{abstract}



\vskip6mm

\section{Introduction}
Let $I=[0,T]$ ($T>0$) and $\mathcal{H}$ a separable Hilbert space.
Consider, for each $(t,x)\in I\times \mathcal{H}$, a maximal
monotone operator of $\mathcal{H}$; $A(t,x)$, for which  the
set-valued map $t \mapsto A(t,x) $ is right-continuous with bounded
variation (BVRC) w.r.t the time $t$, and Lipschitz continuous w.r.t
the state variable $x$,
  in the sense  that
 there exists a function  $\rho:  I \rightarrow [0, \infty[$,
which is right-continuous on $[0, T[$,   and a nonnegative real
number $\gamma$ such that for $0\leq s \leq t \leq T,\;x, y\in
\mathcal{H}$,
$$ dis(A(t,x), A(s,y)) \leq d\rho(]s, t])+ \gamma\|x-y\|= \rho(t)-\rho(s)+\gamma\|x-y\|, $$
where  $dis(\cdot, \cdot)$ is the pseudo-distance between maximal
monotone operators introduced by Vladimirov \cite{Vla}; see relation
\eqref{2.1}, and $f:I\times \mathcal{H}\times
\mathcal{H}\longrightarrow \mathcal{H}$ a   Carath\'eodory mapping
satisfying a linear growth condition.

Let $\lambda$ be the Lebesgue measure on $I$ and $d\rho$  the
Stieljes-measure associated with $\rho$.  We set $\nu := \lambda
+d\rho$ and $\frac{ d\lambda}{d\nu}$ the density of $ \lambda$ w.r.t
 $\nu$.

 In this paper, we are mainly interested by the existence of
  BVRC  mappings $u, v: I\longrightarrow \mathcal{H}$
 satisfying
$$
(P_f)
\begin{cases}
u(0)=u_0, v(0)=v_0\in D(A(0, u_0));\\
 v(t)\in D(A(t, u(t)))\;\;\;\forall t\in I;\\
\displaystyle\frac{du}{d\nu }(t)=v(t)\;\;\;d\nu-a.e.\,t\in
   I;\\
   -\displaystyle\frac{dv}{d\nu }(t)\in A(t,u(t)) v(t)+ f(t,u(t), v(t))  \frac{d\lambda} {d\nu}(t)  \;\;\;d\nu-a.e.\,t\in
   I.
\end{cases}
$$
This work constitutes an extension to the second-order of few
results dealing with first-order BVRC time-dependent maximal
monotone operators differential inclusions, established in \cite{Az,
ACM1, KM}.

For some investigations on first-order differential inclusions
governed by time-dependent and time and state dependent maximal
monotone operators, where the variation of these operators is
absolutely continuous or continuous with bounded variation; we refer
for instance to \cite{ABCM1, ABCM2, AB, ACM, BAC, CIT, Ken, KM, Le,
P, SAM, TOL1, Vla}. Concerning the second order, we refer to the
recent works \cite{CMRF, Saidi1}.

It is well known that the normal cone to closed and convex sets is a
maximal monotone operator, so  differential inclusions governed by
this normal cone, called "sweeping processes" (\cite{M, M1, Mon})
are particular cases of the evolutions with general maximal monotone
operators.  Existence of  BVRC solutions to the first-order sweeping
process
 was   discussed in many interesting papers; we cite
for instance \cite{AHT, ANT, BCG, BCSS, CM, ET1, M2, NNT, Thi3,
TOL2}. However, to the best of our knowledge, there are no results
for the BVRC second-order sweeping process in comparison with the
vast literature concerning absolutely continuous and continuous with
bounded variation solutions; see \cite{AL, AL1, AN, AA1, AA2, AACM,
Az1, AI, Bou, BL, C1, CIY, NNT, N, S}.

The paper is organized as follows. In section 2, we give notations
and recall the needed preliminary results. In section 3, we prove
the existence   of   right-continuous with bounded variation
solutions to the evolution problem $(P_f)$ when $t \mapsto
A(t,\cdot) $ is BVRC and the perturbation $f: I\times
\mathcal{H}\times \mathcal{H}\longrightarrow \mathcal{H}$ is
measurable on $I$ and Lipshitz-continuous on  $\mathcal{H}\times
\mathcal{H}$.

\section{Notations and  Preliminaries} Throughout the paper, we will denote by $\langle \cdot ,\cdot \rangle$ the inner product
of $\mathcal{H}$,   $\|\cdot\|$ the associated norm and
$r\overline{B}_{\mathcal{H}}$ the closed
  ball of center $0$ and radius $r>0$. We will denote by $\mathcal{B}(I)$
the Borel $\sigma$-algebra on $I$. The identity mapping of
$\mathcal{H}$ will be denoted by $Id_{\mathcal{H}}$. For a subset
$K$ of $\mathcal{H}$, $\overline{co}(K)$ will be the closed convex
hull of $K$, which is characterized by:
\begin{equation}\label{co}
\overline{co}(K)=\big\{y\in \mathcal{H}:\;\langle y, x\rangle\leq
\sup_{z\in K}\langle z, x\rangle\;\forall x\in K\big\}.
\end{equation}


If $\mu$ is a positive measure on $I$, we will denote by $L^p(I,
\mathcal{H}; \mu)$ $p \in [1, +\infty[$, the Banach space of classes
of $p$-$\mu$-integrable mappings, equipped with its classical norm
$\|\cdot\|_p$.
\vskip2mm The following results concern  definitions and some
properties of functions with bounded variation and general vector
measures, taken from references \cite{Mon, M2, M3, MV}.

Let $u:I\longrightarrow \mathcal{H}$. The variation of $u$
 in $I$ is the nonnegative extended real number
 \begin{equation}\label{1}
 var(u ; I):=\sup\sum_{k=1}^{n}\|u(t_k)-u(t_{k-1})\|,
 \end{equation}
 where the supremum is taken in $[0,+\infty]$ w.r.t all the finite
 sequences $t_0<t_1<\cdots<t_n$ of points of $I$ ($n$ is arbitrary).
 The function $u$ is said with bounded variation (BV) if and only if
 $var(u ; I)<+\infty$. In this
 case we have
 \begin{equation*}
 \lim_{s\uparrow t}\|u(t)-u(s)\|=\lim_{s\uparrow t} var(u; [s,t]).
 \end{equation*}

 Now, let $\mu$ be a positive Radon measure on $I$ and $\hat{\mu}$
 be an $\mathcal{H}$-valued measure on
 $I$ admitting a density $\frac{d \hat{\mu}}{d\mu}\in
 L^1_{loc}(I,\mathcal{H};\mu)$. Then for $d \mu$-almost every $t\in
 I$, we have
 \begin{equation}\label{2}
\frac{d \hat{\mu}}{d\mu}(t)=\lim_{\varepsilon\downarrow 0}\frac{d
\hat{\mu}([t, t+\varepsilon])}{d\mu([t,
t+\varepsilon])}=\lim_{\varepsilon\downarrow 0}\frac{d
\hat{\mu}([t-\varepsilon, t])}{d\mu([t-\varepsilon, t])}.
 \end{equation}
 The measure $\hat{\mu}$ is absolutely continuous w.r.t $\mu$ if and
 only if $\hat{\mu}=\frac{d \hat{\mu}}{d\mu}\mu$, i.e.,  $\frac{d
 \hat{\mu}}{d\mu}$ is a density of $\hat{\mu}$  w.r.t $\mu$. In this
 case a mapping $u:I\longrightarrow \mathcal{H}$ is
 $\hat{\mu}$-integrable if and only if the mapping $t\mapsto u(t)\frac{d
 \hat{\mu}}{d\mu}(t)$ is $\mu$-integrable, and we have
 \begin{equation}\label{7}
 \int_I u(t) d\hat{\mu}(t)=\int_I u(t)\frac{d
 \hat{\mu}}{d\mu}(t) d\mu(t).
 \end{equation}
If $u:I\longrightarrow \mathcal{H}$ is BVRC and $du$ is its
 differential measure, then we have
 \begin{equation}\label{3}
 u(t)=u(s)+\int_{]s,t]} du\;\;\;\;\;\forall s, t\in I\;(s\leq t).
 \end{equation}
 Conversely, if there exists $v\in L^1(I,\mathcal{H};\mu)$ such that
 $u(t)=u(0)+\int_{]0, t]} v \,d \mu$ for all $t\in I$, then $u$ is
 BVRC and $du=v \, d\mu$, that is $v$ is a density of the vector
 measure $du$ w.r.t the measure $\mu$. So that by \eqref{2}, we get
 \begin{equation}\label{4}
v(t)=\frac{d u}{d\mu}(t)=\lim_{\varepsilon\downarrow 0}\frac{ du
([t, t+\varepsilon])}{d\mu([t,
t+\varepsilon])}=\lim_{\varepsilon\downarrow 0}\frac{ du
([t-\varepsilon, t])}{d\mu([t-\varepsilon, t])}.
 \end{equation}
 If $\mu(\{t\})>0$, this last relation shows that
 \begin{equation}\label{5}
 v(t)=\frac{d u}{d\mu}(t)=\frac{d
 u(\{t\})}{d\mu(\{t\})}\;\;\;\textmd{and}\;\;\;\frac{d
 \lambda}{d\mu}(t)=0.
 \end{equation}

Next, we give definition and some properties of maximal monotone
operators. We refer the reader to \cite{Ba, brezis, Vra} for more
details.

 Let $A: \mathcal{H}
\rightrightarrows \mathcal{H}$ be a set-valued map. We denote by
$D(A)$ and $Gr(A)$ its domain  and graph. We say that the operator
$A$ is monotone, if $\langle y_1 -y_2, x_1 - x_2\rangle\ge 0$ for
all $(x_i, y_i)\in Gr(A)$ $(i=1, 2)$, and
  we say that $A$ is a maximal monotone operator of $\mathcal{H}$,
  if it is monotone and
its graph could not be contained strictly in the graph of any other
monotone operator.

If $A$ is a maximal monotone operator of $\mathcal{H}$, then  for
every $x \in D(A)$, $A(x)$ is nonempty closed and convex. We denote
the projection of the origin on the set $A(x)$ by $A^0(x)$.

 For $\eta > 0$, we denote by $J_{\eta}^A= (Id_{\mathcal{H}} +
\eta A)^{-1}$  the resolvent of  $A$. 
It is
 well-known that
  this operator
is single-valued and defined on all of $\mathcal{H}$, furthermore
$J^A_{\eta}(x)\in D(A)$, for all $x\in \mathcal{H}$.

Let $A: D(A)\subset \mathcal{H}\rightrightarrows \mathcal{H}$ and
$B: D(B)\subset \mathcal{H}\rightrightarrows \mathcal{H}$ be two
maximal monotone operators. Then,  $dis(A, B)$ is the Vladimirov's
pseudo-distance between $A$ and $B$, defined as follows
\begin{equation}\label{2.1} dis(A,B)=\sup\bigg\{\frac{\langle y-y', x'-x\rangle}{1+\|
y\|+\| y'\|}:\;(x,y)\in Gr(A),\;(x',y')\in
Gr(B)\bigg\}.\end{equation}

The following lemmas are  needed for our proof. We refer to
\cite{KM} for their proofs.
\begin{lem}\label{lem2.1}
Let $A$ be a maximal monotone operator of $\mathcal{H}$. If $x\in
\overline{D(A)}$ and  $y\in \mathcal{H}$ are such that
$$ \langle A^0(z)-y, z-x\rangle\geq 0\;\;\forall z\in D(A),$$
then $x\in D(A)$ and $y\in A(x)$.
\end{lem}
\begin{lem}\label{lem2.2}
Let $A_n$ $(n\in \mathbb{N})$ and $A$ be maximal monotone operators
of $\mathcal{H}$ such that $dis(A_n, A)\to 0$. Suppose also that
$x_n\in D(A_n)$ with $x_n\to x$ and $y_n\in A_n(x_n)$ with $y_n\to
y$ weakly for some $x, y\in E$. Then $x\in D(A)$ and $y\in A(x)$.
\end{lem}
\begin{lem}\label{lem2.3}
Let $A$ and $B$ be maximal monotone operators of $\mathcal{H}$. Then\\
1) for $\eta>0$ and $x\in D(A)$
$$ \|x-J_{\eta}^{B}(x)\| \leq \eta\|A^0(x)\|+dis(A,B)+\sqrt{\eta\big(1+\|A^0(x)\|\big)dis(A, B)}.$$
2) For $\eta>0$ and $x, x'\in E$
$$ \|J_{\eta}^{A}(x)-J_{\eta}^{A}(x')\|\leq
\|x-x'\|.$$
\end{lem}
\begin{lem}\label{lem2.4}
Let $A_n$ $(n\in \mathbb{N})$ and $A$ be maximal monotone operators
of $\mathcal{H}$ such that $dis(A_n, A)\to 0$ and $\|A^0_n(x)\|\leq
c(1+\|x\|)$ for some $c>0$, all $n\in \mathbb{N}$ and $x\in D(A_n)$.
Then for every $z\in D(A)$ there exists a sequence $(\zeta_n)$ such
that
 \begin{equation*}\label{2.2}\zeta_n\in
D(A_n),\;\;\;\zeta_n\to z\;\;\textmd{and}\;\;A_n^0(\zeta_n)\to
A^0(z).\end{equation*}
\end{lem}

We end this section by the following discrete Gronwall's lemma.

\begin{lem}\label{lem2.5}
Let $(\alpha_i)$, $(\beta_i)$, $(\gamma_i)$ and $(a_i)$ be sequences
of nonnegative real numbers such that $a_{i+1}\leq
\alpha_i+\beta_i\big(a_0+a_1+\cdots+a_{i-1}\big)+(1+\gamma_i)a_i$
for $i\in \mathbb{N}$. Then
$$ a_j\leq \bigg(a_0+\sum_{k=0}^{j-1}
\alpha_k\bigg)\exp\bigg(\sum_{k=0}^{j-1}
\big(k\beta_k+\gamma_k\big)\bigg)\;\;\textmd{for}\;j\in
\mathbb{N}^*.$$
\end{lem}
\vskip2mm

We will establish our main result under the following hypotheses.\\
$(H_1)$ There exist a function $\rho:I\longrightarrow [0, +\infty[$,
which is right-continuous on $[0, T[$ and nondecreasing with
$\rho(0) = 0$ and $\rho(T)<+\infty$, and nonnegative constant
$\gamma$ such that
\begin{equation*}dis(A(t,x), A(s,y))\leq d\rho(]s,
t])+\gamma\|x-y\|\frac{d\lambda}{d\nu}(s)\;\;\textmd{for}\;\;0\leq
s\leq t\leq T,\;x,y\in \mathcal{H},\end{equation*}  and since
$\frac{d\lambda}{d\nu}(s)\leq 1$, $(H_1)$ implies that
\begin{equation*}dis(A(t,x), A(s,y))\leq d\rho(]s,
t])+\gamma\|x-y\|\;\;\textmd{for}\;\;0\leq s\leq t\leq T,\;x,y\in
\mathcal{H}.\end{equation*}
 $(H_2)$ There
exists a nonnegative real constant $c$ such that
\begin{equation*}\|A^0(t,x)y\|\leq c(1+\| x\|+\| y\|)\;\;\textmd{for}\;\;t\in
I,\;x\in\mathcal{H},\;y\in D(A(t,x)).\end{equation*}
 $(H_3)$ For any bounded subset $E$ of $\mathcal{H}$, $D(A(I\times E))$
 is relatively ball-compact, i.e., its intersection with any closed ball of
 $\mathcal{H}$ is relatively compact.\\
$(H_4)$ $t\mapsto f(t,\cdot,\cdot)$ is $\mathcal{B}(I)$-measurable.\\
$(H_5)$ There  is a nonnegative function $k\in L^1(I,\mathbb{R};
\lambda)$ such that
$$ \|f(t,x,y) - f(t,x',y')\|\leq k(t)(\|x-x'\|+\|y-y'\|)\;\;\;\forall t\in I, \;\forall (x,y), (x',y')\in
\mathcal{H}\times \mathcal{H}.$$
 $(H_6)$ There exists a nonnegative real constant $m$
such that
$$ \|f(t,x,y)\|\leq m(1+\|x\|
+\|y\|)\;\;\;\forall (t,x,y)\in I\times \mathcal{H}\times
\mathcal{H}.$$

\section{Main result}
Now we are able to state our main theorem.
\begin{thm}\label{Theorem 3.1}    Let for
every $(t,x)\in I\times \mathcal{H}$, $A(t,x):D(A(t,x))\subset
\mathcal{H}\rightrightarrows \mathcal{H}$ be a maximal monotone
operator satisfying $(H_1)$, $(H_2)$ and $(H_3)$. Let $f:I\times
\mathcal{H}\times \mathcal{H}\longrightarrow \mathcal{H}$  such that
$(H_4)$, $(H_5)$ and $(H_6)$ are satisfied. Then for any
$(u_0,v_0)\in \mathcal{H}\times D(A(0, u_0))$, there exists a BVRC
solution $ (u, v) :I \longrightarrow \mathcal{H}\times \mathcal{H}$
to  problem $(P_f)$.
\end{thm}

\proof
 Following
Castaing et al \cite{CM}, we can choose a sequence
$(\varepsilon_n)_{n\geq 1}\subset ]0, 1]$ such that $\varepsilon_n
\downarrow 0$ and a partition $0=t_0^n<t_1^n<\cdots<t_{q_n}^n=T$ of
$I$ (with $\underset{n\to\infty}{\lim}q_n=+\infty$), for which we
have
\begin{equation}\label{3.1}
\nu(]t^n_i, t^n _{i+1}])=|t_{i+1}^n-t_i^n|+d\rho(]t_i^n,
t_{i+1}^n])\leq
\varepsilon_n\;\;\;\textmd{for}\;i=0,\cdots,q_n-1,\end{equation} we
can take, for example, $\varepsilon_n=\frac{T}{q_n}$. For each $i\in
\{0,\cdots,q_n-1\}$, put $I_i^n=]t_i^n, t_{i+1}^n]$,
\begin{equation}\label{3.2} \delta_{i+1}^n=d\rho(]t_i^n,
t_{i+1}^n]), \;\;\;\; \eta^n_{i+1} = t ^n_{i+1}-t ^n_i, \;\;\;\;
\beta ^n_{i+1}  = \nu(]t^n_i, t^n _{i+1}]).\end{equation} For every
$n\in \mathbb{N}$, put $u_0^n=u_0$, $v_0^n=v_0\in D(A(0,u_0))$, and
let us define by induction,  sequences $(v_i^n)_{0\leq i\leq
q_n-1}$, $(u_i^n)_{0\leq i\leq q_n-1}$ such that
\begin{equation}\label{3.3} v_{i+1}^n= J^n_{i+1}
\Big(v_i^n-\int_{t_i^n}^{t_{i+1}^n} f(s, u_i^n,
v_i^n)d\lambda(s)\Big)\end{equation}  where
$J^n_{i+1}:=J^{A(t_{i+1}^n,u_i^n )}_{\beta^n_{i+1}}
=\big(Id_{\mathcal{H}}+ \beta^n_{i+1} A(t_{i+1}^n,
u_i^n)\big)^{-1}$, and
\begin{equation}\label{3.2'}
u_{i+1}^n=u_i^n+\beta^n_{i+1}v_{i+1}^n.
\end{equation}
  From the properties of the resolvent,  we have that  $v^n_{i+1}
\in D(A(t_{i+1}^n, u_i^n))$ and
\begin{equation}\label{3.5} -\frac{1} {\beta^n_{i+1}}\Big(v^n_{i+1}-
v_i^n+\int_{t_i^n}^{t_{i+1}^n} f(s, u_i^n, v_i^n)d\lambda(s)\Big)
\in A( t_{i+1}^n, u_i^n) v^n _{i+1}. \end{equation} For $t\in
[t^n_i, t ^n_{i+1}[$,  $i = 0,\cdots,q_n-1$, set
\begin{equation}\label{3.7}
v_n(t)=   v_i^n + \frac{ \nu(]t^n_i, t])  }{ \nu(]t^n_i, t
^n_{i+1}])} \Big(v_{i+1}^n-v_i^n+\int_{t_i^n}^{t_{i+1}^n} f(s,u_i^n,
v_i^n) d\lambda(s)\Big)  - \int_{t_i^n}^{t}
f(s,u_i^n,v_i^n)d\lambda(s),\end{equation}
\begin{equation}\label{3.7'}
u_n(t)=   u_i^n +  \frac{ \nu(]t^n_i, t])  }{ \nu(]t^n_i, t
^n_{i+1}])}  (u_{i+1}^n-u_i^n),
\end{equation} and
$v_n(T)=v_{q_n}^n$, $u_n(T)=u_{q_n}^n$. So that $v_n(t_i^n)=v_i^n$,
$u_n(t_i^n)=u_i^n$ and $v_n$, $u_n$ are BVRC mappings on $I$.

\vskip2mm{\bf Step 1.}
 We prove in this step  that the sequences $(v_n)$, $(u_n)$ are bounded in  norm and
 in variation.\\Using Lemma \ref{lem2.3}, we have from
\eqref{3.3}, $(H_1)$,  $(H_2)$ and $(H_6)$, for $i=0,1,\cdots,
q_n-1$ (keeping in mind that $\sqrt{a\,b}\leq \frac{a}{2}+b$ for
nonnegative real numbers $a$, $b$)
\begin{eqnarray*}
\| v_{i+1}^n-v_i^n\| &\leq& \Big\|J_{i+1}^n\Big(v_i^n-
\int_{t_i^n}^{t_{i+1}^n}
f(s,u_i^n,v_i^n)d\lambda(s)\Big)-J_{i+1}^n(v_i^n)\Big\|+\big\|J_{i+1}^n(v_i^n)-v_i^n\big\|\\
&\leq& \int_{t_i^n}^{t_{i+1}^n} \|f(s,u_i^n,v_i^n)\|
d\lambda(s)+\beta_{i+1}^n\|A^0(t_i^n,
u_{i-1}^n)v_i^n\|+dis\big(A(t_{i+1}^n, u_i^n),A(t_i^n,
u_{i-1}^n)\big)\\&+&\sqrt{\beta_{i+1}^n\big(1+\|A^0(t_i^n,
u_{i-1}^n)v_i^n\|\big)dis\big(A(t_{i+1}^n, u_{i}^n),A(t_i^n, u_{i-1}^n)\big)}\\
&\leq& m(1+\|u_i^n\|+\|v_i^n\|)\beta_{i+1}^n
+\big(1+c(1+\|u_{i-1}^n\|+\|v_i^n\|)\big)\beta_{i+1}^n+\gamma\|u_{i}^n-u_{i-1}^n\|\\&+&
\sqrt{\beta_{i+1}^n\big(1+c(1+\|u_{i-1}^n\|+\|v_i^n\|)\big)(\beta_{i+1}^n+\gamma\|u_{i}^n-u_{i-1}^n\|)}\\
&\leq& m(1+\|u_i^n\|+\|v_i^n\|)\beta_{i+1}^n
+\big(1+c(1+\|u_{i-1}^n\|+\|v_i^n\|)\big)\beta_{i+1}^n+\gamma\|u_{i}^n-u_{i-1}^n\|\\
&+&\beta_{i+1}^n\frac{1}{2}\big(1+c(1+\|u_{i-1}^n\|+\|v_i^n\|)\big)+\beta_{i+1}^n+\gamma\|u_{i}^n-u_{i-1}^n\|,
\end{eqnarray*}
that is, \begin{eqnarray}\label{3.8}\| v_{i+1}^n-v_i^n\|&\leq&
(m+1)(1+\|u_i^n\|+\|v_i^n\|)\beta_{i+1}^n+\frac{3}{2}(1+c)(1+\|u_{i-1}^n\|+\|v_{i}^n\|)\beta_{i+1}^n+2\gamma\|u_{i}^n-u_{i-1}^n\|\nonumber\\&=:&
m_1(1+\|u_i^n\|+\|v_i^n\|)\beta_{i+1}^n+c_1(1+\|u_{i-1}^n\|+\|v_{i}^n\|)\beta_{i+1}^n+\gamma_1\|u_{i}^n-u_{i-1}^n\|,\end{eqnarray}
which entails \begin{equation} \label{eq0}\|v_{i+1}^n\|\leq
\|v_i^n\|+m_1(1+\|u_i^n\|+\|v_i^n\|)\beta_{i+1}^n+c_1(1+\|u_{i-1}^n\|+\|v_{i}^n\|)\beta_{i+1}^n+\gamma_1\|u_{i}^n-u_{i-1}^n\|
.\end{equation} On the other hand, we have from \eqref{3.2'},
\begin{eqnarray}\label{eq1}
\|u_{i}^n-u_{i-1}^n\|=\beta_i^n\|v_i^n\|
\end{eqnarray}
and (using \eqref{3.1})
\begin{eqnarray}\label{eq2}
\|u_{i}^n\|&=&\|u_{i-1}^n+\beta_i^nv_i^n\|=\|u_0+\beta_1^nv_1^n+\beta_2^nv_2^n+\cdots+\beta_i^nv_i^n\|\nonumber\\&\leq&
\|u_0\|+\varepsilon_n(\|v_1^n\|+\|v_2^n\|+\cdots+\|v_i^n\|).
\end{eqnarray}
Replacing  these two last relations in \eqref{eq0}, we get by the
use of \eqref{3.1},
\begin{eqnarray*}
\|v_{i+1}^n\|&\leq&\|v_i^n\|+m_1\Big(1+\|u_0\|+\varepsilon_n(\|v_1^n\|+\|v_2^n\|+\cdots+\|v_i^n\|)+\|v_i^n\|\Big)\varepsilon_n\\&+&
c_1\big(1+\|u_0\|+\varepsilon_n(\|v_1^n\|+\|v_2^n\|+\cdots+\|v_{i-1}^n\|)+\|v_i^n\|\big)\varepsilon_n+\gamma_1\|v_i^n\|\varepsilon_n\\&\leq&
(m_1+c_1)(1+\|u_0\|)\varepsilon_n+(m_1+c_1)(\|v_1^n\|+\|v_2^n\|+\cdots+\|v_{i-1}^n\|)\varepsilon^2_n
\\&+&\Big(1+(m_1\varepsilon_n+m_1+c_1+\gamma_1)\varepsilon_n\Big)\|v_i^n\|.
\end{eqnarray*}
Applying Lemma \ref{lem2.5}, we get for all $i\leq q_n$,
\begin{eqnarray*}\|v_i^n\|&\leq&
\Big(\|v_0\|+\sum_{j=0}^{i-1}
(m_1+c_1)(1+\|u_0\|)\varepsilon_n\Big)\exp\Big(\sum_{j=0}^{i-1}
\big(j(m_1+c_1)\varepsilon^2_n+(m_1\varepsilon_n+m_1+c_1+\gamma_1)\varepsilon_n\big)\Big)\\
&\leq&\Big(\|v_0\|+(m_1+c_1)(1+\|u_0\|)T)\Big)\exp\Big((2+T)m_1+2c_1+\gamma_1
\Big)=:m_2,
\end{eqnarray*}
using this last estimate in \eqref{eq2}, we obtain
$$\|u_i^n\|\leq \|u_0\|+m_2 T,$$ and replacing in
 \eqref{eq1} and in \eqref{3.8}, it comes that
 $$\|u_{i+1}^n-u_i^n\|\leq m_2\beta_{i+1}^n,$$ and
$$\|v_{i+1}^n-v_i^n\|\leq(m_1+c_1) (1+\|u_0\|+m_2T+m_2)\beta_{i+1}^n+\gamma_1m_2\beta_{i}^n=:m_3\beta_{i+1}^n+m_4\beta_{i}^n.$$
Consequently,
\begin{equation}\label{3.9}\|u_i^n\|\leq M,\;\;\;\|v_i^n\|\leq M\;\;\forall\, 0\leq i \leq
q_n\end{equation} and for all $0\leq i<q_n-1$,
 \begin{equation}\label{3.9'}\|u_{i+1}^n-u_i^n\|\leq M\nu(]t_i^n,
t_{i+1}^n]),\;\;\;\|v_{i+1}^n-v_i^n\|\leq M\big(\nu(]t_i^n,
t_{i+1}^n])+\nu(]t_{i-1}^n, t_{i}^n])\big),
\end{equation}
where $M=\max(m_2, m_3, m_4)$. Now, for $t\in [t_i^n, t_{i+1}^n[$,
we have from \eqref{3.1}, \eqref{3.7}, \eqref{3.9}, \eqref{3.9'} and
$(H_6)$,
\begin{equation*}\label{3.10}
\|v_n(t)-v_i^n\| \leq
2M\varepsilon_n+2m(1+2M)\varepsilon_n\leq2M+2m(1+2M)=:M_1,\end{equation*}
and from \eqref{3.1}, \eqref{3.7'},  \eqref{3.9} and \eqref{3.9'},
we get
\begin{equation*}\label{3.10'}
\|u_n(t)-u_i^n\| \leq M.
\end{equation*}
These last relations give us, for all $n\in \mathbb{N}$,
\begin{equation*}\|v_n(t)\|\leq M+M_1=:M_2\;\;\;\textmd{and}\;\;\;\|u_n(t)\|\leq 2M\;\;\;\forall t\in I,\end{equation*}  that is
\begin{equation*}\label{3.11}
\sup_n\|v_n\|:=\sup_n\big(\sup_{t\in I}\|v_n(t)\|\big)\leq
M_2\;\;\;\textmd{and}\;\;\;\sup_n\|u_n\|\leq 2M.
\end{equation*}
That is $(u_n)$, $(v_n)$ are bounded in norm and in variation.
\vskip2mm

{\bf Step 2.} Convergence of the sequences $(u_n)$, $(v_n)$ and $(\frac{dv_n}{d\nu})$.\\
Define the functions $\phi_n,\;\theta_n: I\longrightarrow I$ by
$$\phi_n(t)=t_{i}^n,\;\;\;\theta_n(t)=t_{i+1}^n\;\;\textmd{for}\;\; t\in ]t_i^n,
t_{i+1}^n],\;i=0,1,\cdots,q_n-1,\;\; \textmd{and}\;\;
\phi_n(0)=\theta_n(0)=0,
$$ and set for all $t\in I$,
\begin{equation*}
B_n(t)=\sum_{i=0}^{q_n-1}\frac{1}{\nu(]t_i^n, t_{i+1}^n])}\Big(
v_{i+1}^n-v_i^n+\int_{t_i^n}^{t_{i+1}^n} f(s,u_i^n,
v_i^n)d\lambda(s)\Big)\chi_{]t_i^n, t_{i+1}^n]}(t),
\end{equation*}
where $\chi_{J}(\cdot)$ is the characteristic function of the set
$J\subset I$. Whence, by \eqref{3.7}, it is clear that
\begin{equation*}
v_n(t)=v_0+\int_{]0,t]} B_n(\tau)d\nu(\tau)-\int_{]0,t]} f(\tau,
u_n(\phi_n(\tau)), v_n(\phi_n(\tau)))d\lambda(\tau).
\end{equation*}
Since $\frac{d\lambda}{d\nu}$ is a density of $\lambda$ w.r.t $\nu$,
 then by \eqref{7}, we get for every $t\in I$,
\begin{equation*}
v_n(t)=v_0+\int_{]0,t]} \Big(B_n(\tau)-f(\tau, u_n(\phi_n(\tau)),
v_n(\phi_n(\tau)))\frac{d\lambda}{d\nu}(\tau)\Big)d\nu(\tau),
\end{equation*}
so that, $\frac{dv_n}{d\nu}$ is a density of the vector measure
$dv_n$ w.r.t $\nu$ and for $\nu$-almost every $t\in I$,
\begin{equation*}\label{8}\frac{dv_n}{d\nu}(t)=B_n(t)-f(t,
u_n(\phi_n(t)),
v_n(\phi_n(t)))\frac{d\lambda}{d\nu}(t),\end{equation*} that is, by
the definition of $B_n$, for $\nu$-almost every $t\in I$,
\begin{eqnarray}\label{8}
&&\frac{dv_n}{d\nu}(t)+f(t,
u_n(\phi_n(t)),v_n(\phi_n(t)))\frac{d\lambda}{d\nu}(t)\nonumber\\&=&\sum_{i=0}^{q_n-1}\frac{1}{\nu(]t_i^n,
t_{i+1}^n])}\Big( v_{i+1}^n-v_i^n+\int_{t_i^n}^{t_{i+1}^n}
f(s,u_i^n, v_i^n)d\lambda(s)\Big)\chi_{]t_i^n, t_{i+1}^n]}(t).
\end{eqnarray}
 From $(H_6)$, \eqref{3.9} and \eqref{8}, it results that
 \begin{equation}\label{9}
 \Big\|\frac{dv_n}{d\nu}(t)+f(t,
u_n(\phi_n(t)), v_n(\phi_n(t)))\frac{d\lambda}{d\nu}(t)\Big\|\leq
M+m(1+2M)=:M_3\;\;\;\nu-a.e.\,t\in I.
 \end{equation}
On the other hand, we know that $\nu=d\rho+\lambda$, whence
$\frac{d\lambda}{d\nu}(t)\leq 1$. So that, by $(H_6)$ and
\eqref{3.9}, it follows that
\begin{equation}\label{10}
 \Big\|f(t,
u_n(\phi_n(t)), v_n(\phi_n(t)))\frac{d\lambda}{d\nu}(t)\Big\|\leq
m(1+2M)=:M_4\;\;\;\nu-a.e.\,t\in I,
 \end{equation}
 consequently, there is a Borel subset $J\subset I$ with $\nu(J)=0$,
 such that
 \begin{equation}\label{11}
 \Big\|\frac{dv_n}{d\nu}(t)\Big\|\leq M_3+M_4=:M_5\;\;\;\forall t\in I\setminus J.
 \end{equation}
 Now, observe that by \eqref{3.5} and \eqref{8}, for each $n\in \mathbb{N}$, there is a
 Borel subset $J_n\subset I$ with $\nu(J_n)=0$ such that
 \begin{equation}\label{12}
 \frac{dv_n}{d\nu}(t)+f(t,
u_n(\phi_n(t)),v_n(\phi_n(t)))\frac{d\lambda}{d\nu}(t)\in
-A(\theta_n(t), u_n(\phi_n(t)))v_n(\theta_n(t))\;\;\;\forall t\in
I\setminus J_n,
\end{equation}
moreover,
\begin{equation}\label{13}
v_n(\theta_n(t))\in
D\big(A(\theta_n(t),u_n(\phi_n(t)))\big)\;\;\;\forall t\in I.
\end{equation}

We will show in the following, that $(v_n)$ converges uniformly on $I$.\\
On the one hand, since $\frac{dv_n}{d\nu}$ is a density of the
vector measure $dv_n$ w.r.t the measure $\nu$, we have from
\eqref{3}, for all $s, t\in I$ $(s\leq t)$,
\begin{eqnarray*}
v_n(t)=v_n(s)+\int_{]s, t]} \frac{dv_n}{d\nu}(s)\,d\nu(s),
\end{eqnarray*}
that is, using relation \eqref{11}, we get
\begin{equation}\label{eq3}\|v_n(t)-v_n(s)\|\leq M_5 \nu(]s,
t]),\end{equation} which means that $(v_n)$ is equi-right-continuous
with bounded variation. On the other hand, from \eqref{3.9} and
\eqref{13}, it is clear that for any $t\in I$, we have
$$\big(v_n(\theta_n(t))\big)\subset D(I\times
M\overline{B}_{\mathcal{H}})\cap M \overline{B}_{\mathcal{H}},$$
which shows that the sequence $\big(v_n(\theta_n(t))\big)$ is
relatively compact (see $(H_3))$. From \eqref{eq3}, we have that for
all $t\in I$,
\begin{equation}\label{eq4}
\|v_n(\theta_n(t))-v_n(t)\|\leq M_5d\nu(]t, \theta_n(t)])\to
0\;\;\textmd{as}\;\;n\to\infty,
\end{equation}
whence, we conclude that $(v_n(t))$ is also relatively compact.
Using Helly-Banach's Theorem (see \cite{Por}), we may assume the
existence of a BVRC mapping $v$ such that $(v_n(t))$ converges
strongly to $v(t)$ for every $t\in I$. This together with
\eqref{eq3}, show that
\begin{equation}\label{eq3'}\|v(t)-v(s)\|\leq M_5 \nu(]s,
t]),\end{equation} and from \eqref{eq4}
\begin{equation}\label{eq4'}
\|v_n(\theta_n(t))-v(t)\|\leq
\|v_n(\theta_n(t))-v_n(t)\|+\|v_n(t)-v(t)\|\to
0\;\;\textmd{as}\;\;n\to\infty.
\end{equation}
As regards the sequence $(u_n)$, let us set for all $t\in I$,
\begin{equation*}
\alpha_n(t)=\sum_{i=0}^{q_n-1}\frac{1}{\nu(]t_i^n, t_{i+1}^n])}
(u_{i+1}^n-u_i^n)\chi_{]t_i^n, t_{i+1}^n]}(t),
\end{equation*}
and remark that by relation \eqref{3.2'},
$\alpha_n(t)=v_n(\theta_n(t))$ for $\nu$-a.e.  $t\in I$.
Consequently, from \eqref{3.7'}, it comes that for all $t\in I$,
\begin{equation}\label{eq5}
u_n(t)=u_0+\int_{]0, t]}\alpha_n(\tau) d\nu(\tau)=u_0+\int_{]0,
t]}v_n(\theta_n(\tau))d\nu(\tau).
\end{equation}
Using relations \eqref{3.9} and \eqref{eq4} and the strong
convergence in $\mathcal{H}$ of $(v_n(t))$ to $v(t)$, we obtain by
Lebesgue dominated convergence, that
$$ \lim_{n\to\infty} u_n(t)=u_0+\int_{]0,
t]}v(\tau)d\nu(\tau)=:u(t)$$ so that, $u$ is a BVRC mapping,
$(u_n(t))$ converges strongly to  $u(t)$ for all $t\in I$, and
$\frac{du}{d\nu}$ is a density of the vector measure $du$ w.r.t
$\nu$ and for $\nu$-almost every $t\in I$,
\begin{equation}\label{noter}
\frac{du}{d\nu}(t)=v(t).
\end{equation}

Now, observe that the sequence $(\frac{dv_n}{d\nu})$ is bounded in
$L^2(I, \mathcal{H}; \nu)$ due to relation \eqref{11}, so that it
converges weakly in $L^2(I, \mathcal{H}; \nu)$ to some mapping $w\in
L^2(I, \mathcal{H}; \nu)$. In particular, for all $t\in I$
\begin{equation*}
\int_{]0, t]} \frac{dv_n}{d\nu}(s)\, d\nu(s)\longrightarrow\int_{]0,
t]} w(s)\, d\nu(s)\;\;\;\textmd{weakly in }\;\mathcal{H}.
\end{equation*}
Since $\frac{dv_n}{d\nu}$ is a density of the vector measure $dv_n$
w.r.t the measure $\nu$, we have for all $t\in I$,
\begin{eqnarray*}
v_n(t)=v_0+\int_{]0, t]} \frac{dv_n}{d\nu}(s)\,d\nu(s)
\end{eqnarray*}
and since for any $t\in I$, $(v_n(t))$ converges strongly in
$\mathcal{H}$ and then weakly to $v(t)$, we deduce by what precedes
that
\begin{equation*}
v(t)=v_0+\int_{]0, t]} w(s)\,d\nu(s)\;\;\;\forall t\in I.
\end{equation*}
This shows that $\frac{dv}{d\nu}=w$ for $\nu$-a.e. $t\in I$.
Consequently, $(\frac{dv_n}{d\nu})$ converges weakly in $L^2(I,
\mathcal{H}; \nu)$ to $\frac{dv}{d\nu}$.

\vskip2mm

{\bf Step 3.} Existence of solution.\\
First, observe that from \eqref{13},  $ v_n(\theta_n(t))\in
D\big(A(\theta_n(t),u_n(\phi_n(t)))\big)$ for all $t\in I $. Also,
from $(H_1)$, we have that for all $t\in I$,
$$dis\big(A(\theta_n(t), u_n(\phi_n(t))), A(t, u(t))\big)\leq
|\rho(\theta_n(t))-\rho(t)|+\gamma\|u_n(\phi_n(t))-u(t)\|\frac{d\lambda}{d\nu}(t),$$
but from \eqref{eq5}, it is clear that
$$\|u_n(\phi_n(t))-u(t)\|\leq \|u_n(\phi_n(t))-u_n(t)\|+\|u_n(t)-u(t)\|\leq M\nu(]\phi_n(t), t])+\|u_n(t)-u(t)\|,$$
so, from \eqref{5}, it comes that
\begin{equation}\label{eq5'}\|u_n(\phi_n(t))-u(t)\|\frac{d\lambda}{d\nu}(t)\leq
\big(M\nu(]\phi_n(t),
t])+\|u_n(t)-u(t)\|\big)\frac{d\lambda}{d\nu}(t)\longrightarrow
0\;\textmd{as}\;n\to\infty.\end{equation} This shows that
\begin{equation}\label{eq6}
dis\big(A(\theta_n(t), u_n(\phi_n(t))), A(t,
u(t))\big)\longrightarrow 0\;\textmd{as}\;n\to\infty.
\end{equation}
On the other hand, by \eqref{3.9} and $(H_2)$, we have that the
sequence $\big(A^0\big(\theta_n(t),
u_n(\phi_n(t))\big)v_n(\theta_n(t))\big)$ is bounded in
$\mathcal{H}$, and hence it is weakly relatively compact. Using this
fact, \eqref{eq4'} and \eqref{eq5'}, we conclude by Lemma
\ref{lem2.2}, that $v(t)\in D(A(t, u(t)))$ for all $t\in I$.

Next, remember that the sequence $(\frac{dv_n}{d\nu})$ converges
weakly in $L^2(I, \mathcal{H}; \nu)$ to $\frac{dv}{d\nu}$, so that
by Mazur's theorem, there is a sequence $(\xi_n)$ such that for each
$n\in\mathbb{N}$, $\xi_n\in co\{\frac{dv_k}{d\nu};\;k\geq n\}$, and
$(\xi_n)$ converges strongly in $L^2(I, \mathcal{H}; \nu)$ to
$\frac{dv}{d\nu}$. Whence, there is a subsequence $(\xi_{n_j})$,
which converges $\nu$-almost every where to $\frac{dv}{d\nu}$. This
means the existence of a Borel subset $J'\subset I$, with
$\nu(J')=0$ and for $t\in I\setminus J'$,
$$\xi_{n_j}(t)\longrightarrow \frac{dv}{d\nu}(t)\in\bigcap_{n}
\overline{co} \big\{\frac{dv_{k}}{d\nu};\;k\geq n\big\}.$$ This
implies, by \eqref{co}, that for any fixed $\eta\in \mathcal{H}$,
\begin{equation}\label{17}
\big\langle \frac{dv}{d\nu}(t), \eta
\big\rangle\leq\limsup_{n\to\infty}\big\langle \frac{dv_n}{d\nu}(t),
\eta \big\rangle.
\end{equation}
On the other hand, we have from $(H_5)$ and \eqref{eq3}
\begin{eqnarray*}
&&\big\|f\big(t,u_n(\phi_n(t)),v_n(\phi_n(t))\big)\frac{d\lambda}{d\nu}(t)-f(t,u(t),v(t))\frac{d\lambda}{d\nu}(t)\big\|
\\&\leq&k(t)\big(\|u_n(\phi_n(t))-u(t)\|+\|v_n(\phi_n(t))-v(t)\|\big)\frac{d\lambda}{d\nu}(t)\nonumber\\
&\leq&k(t)\big(\|u_n(\phi_n(t))-u(t)\|+\|v_n(\phi_n(t))-v_n(t)\|+\|v_n(t)-v(t)\|\big)\frac{d\lambda}{d\nu}(t)\nonumber\\
&\leq&k(t)\big(\|u_n(\phi_n(t))-u(t)\|+M_5\nu(]\phi_n(t),
t])+\|v_n(t)-v(t)\|\big)\frac{d\lambda}{d\nu}(t),
\end{eqnarray*}
whence by the  convergence of $(v_n(t))$ to $v(t)$ and relations
 \eqref{5}, \eqref{eq5'}, it results that for all $t\in I$,
\begin{equation}\label{lip}
f\big(t,u_n(\phi_n(t)),
v_n(\phi_n(t))\big)\frac{d\lambda}{d\nu}(t)\longrightarrow f(t,u(t),
v(t))\frac{d\lambda}{d\nu}(t)\;\;\;\textmd{as}\;\;n\to\infty.
\end{equation}

Now, to prove that $(u, v)$ is a solution to our considered problem,
we will use Lemma \ref{lem2.1}. Since $v(t)\in D(A(t,u(t)))$ for all
$t\in I$, we have to show that for $\nu$-almost every fixed $t\in I$
and for any $z\in D(A(t, u(t)))$
\begin{equation*}
\big\langle
A^0(t,u(t))z+\frac{dv}{d\nu}(t)+f(t,u(t),v(t))\frac{d\lambda}{d\nu}(t),z-v(t)\big\rangle\geq
0.
\end{equation*}
Indeed, let $t\in I$. By $(H_2)$, using Lemma \ref{lem2.4}, we can
ensure the existence of a sequence $(\zeta_n)_n$, such that
\begin{equation}\label{18}
\zeta_n\in
D\big(A(\theta_n(t),u_n(\phi_n(t)))\big),\;\;\;\zeta_n\longrightarrow
z\;\;\textmd{and}\;\; A^0\big(\theta_n(t),u_n(\phi_n(t))\big)
\zeta_n\longrightarrow A^0(t,u(t))z.
\end{equation}
Since $A(\theta_n(t), u_n(\phi_n(t)))$ is monotone, using
\eqref{12}, we have for each $n\in \mathbb{N}$ and $t\in I\setminus
J_n$,
\begin{equation}\label{19}
\big\langle v_n(\theta_n(t))-\zeta_n,
\frac{dv_n}{d\nu}(t)+f(t,u_n(\phi_n(t)),v_n(\phi_n(t)))\frac{d\lambda}{d\nu}(t)+A^0(\theta_n(t),
u_n(\phi_n(t)))\zeta_n\big\rangle\leq 0.
\end{equation}
Next,  let $t\in I\setminus \big(J\cup J'\cup(\underset{n}{\bigcup}
J_n)\big)$. From \eqref{17}, we have
  \begin{equation*}
  \big\langle \frac{dv}{d\nu}(t), v(t)-z\big\rangle\leq\limsup_{n\to\infty}\big\langle \frac{dv_n}{d\nu}(t),
 v(t)-z\big\rangle,
  \end{equation*}
and since from  \eqref{11} and \eqref{19}, we have for each $n\in
\mathbb{N}$,
\begin{eqnarray*}
&&\Big\langle \frac{dv_n}{d\nu}(t),
  v(t)-z\Big\rangle\leq\Big\langle
  \frac{dv_n}{d\nu}(t),v(t)-v_n(\theta_n(t))\Big\rangle+ \Big\langle\frac{dv_n}{d\nu}(t), v_n(\theta_n(t))-\zeta_n\Big\rangle+\Big\langle
  \frac{dv_n}{d\nu}(t),\zeta_n-z\Big\rangle\\&\leq&
  M_5\big(\|v(t)-v_n(\theta_n(t))\|+\|\zeta_n-z\|\big)+\Big\langle
  \frac{dv_n}{d\nu}(t),v_n(\theta_n(t))-\zeta_n\Big\rangle
  \\&\leq&
  M_5\big(\|v(t)-v_n(\theta_n(t))\|+\|\zeta_n-z\|\big)\\&+&\Big\langle f\big(t,u_n(\phi_n(t)), v_n(\phi_n(t))\big)\frac{d\lambda}{d\nu}(t)+A^0(\theta_n(t),
u_n(\phi_n(t)))\zeta_n,\zeta_n-v_n(\theta_n(t))\Big\rangle
\end{eqnarray*}
we conclude, using \eqref{eq4'}, \eqref{lip} and \eqref{18}
that\begin{equation} \big\langle \frac{dv}{d\nu}(t),
v(t)-z\big\rangle\leq\big\langle
f(t,u(t),v(t))\frac{d\lambda}{d\nu}(t)+A^0(t,
u(t))z,z-v(t)\big\rangle,
\end{equation}
which, with the fact that $v(t)\in D(A(t,u(t)))$, $u(0)=u_0$,
$v(0)=v_0$ and \eqref{noter}, means that $(u,v)$ is a BVRC solution
to our problem $(P_f)$. $\hfill\square$ \vskip2mm

\begin{rem}
If the operators depend only on the time and satisfy the following
hypotheses:\\
$(H'_1)$ There exists a function $\rho:I\longrightarrow [0,
+\infty[$, which is right-continuous on $[0, T[$ and nondecreasing
with $\rho(0) = 0$ and $\rho(T)<+\infty$ such that
\begin{equation*}dis(A(t), A(s))\leq d\rho(]s,
t])\;\;\textmd{for}\;\;0\leq s\leq t\leq T.\end{equation*}
 $(H'_2)$ There
exists a nonnegative real constant $c$ such that
\begin{equation*}\|A^0(t)x\|\leq c(1+\| x\|)\;\;\textmd{for}\;\;t\in
I,\;x\in D(A(t)).\end{equation*}
 $(H'_3)$ For any $t\in I$, $D(A(t))$
 is relatively ball-compact.\\
 Our theorem reads as
\end{rem}

\begin{coro}\label{Theorem 3.2}    Let for
every $t\in I$, $A(t):D(A(t))\subset \mathcal{H}\rightrightarrows
\mathcal{H}$ be a maximal monotone operator satisfying $(H'_1)$,
$(H'_2)$ and $(H'_3)$. Let $f:I\times \mathcal{H}\times
\mathcal{H}\longrightarrow \mathcal{H}$  such that $(H_4)$, $(H_5)$
and $(H_6)$ are satisfied. Then for any $(u_0,v_0)\in
\mathcal{H}\times D(A(0))$, there exists a BVRC solution $ (u, v) :I
\longrightarrow \mathcal{H}\times \mathcal{H}$ to  the problem
$$
\begin{cases}
u(0)=u_0, v(0)=v_0\in D(A(0));\\
 v(t)\in D(A(t))\;\;\;\forall t\in I;\\
\displaystyle\frac{du}{d\nu }(t)=v(t)\;\;\;d\nu-a.e.\,t\in
   I;\\
   -\displaystyle\frac{dv}{d\nu }(t)\in A(t) v(t)+ f(t,u(t), v(t))  \frac{d\lambda} {d\nu}(t)  \;\;\;d\nu-a.e.\,t\in
   I.
\end{cases}
$$
\end{coro}
\begin{rem}
In particular, if we consider $A(t,x)=N_{C(t,x)}$, the normal cone,
in the sense of convex analysis, to the closed convex set $C(t,x)$,
then under the following hypotheses:\\ $(H''_1)$ There exists a
function $\rho:I\longrightarrow [0, +\infty[$, which is
right-continuous on $[0, T[$ and nondecreasing with $\rho(0) = 0$
and $\rho(T)<+\infty$ and a nonnegative real constant $\gamma$ such
that
\begin{equation*}
|d_{C(t,x)}(z)-d_{C(s,y)}(z)|\leq d\rho(]s,
t])+\gamma\|x-y\|\frac{d\lambda}{d\nu}(s)\;\;\textmd{for}\;\;0\leq
s\leq t\leq T,\;x,y,z\in \mathcal{H}.
\end{equation*}
$(H''_2)$ For any bounded subset $E$ of $\mathcal{H}$, $C(I\times
E)$
 is relatively ball-compact.\\
 We can derive the following corollary.
\end{rem}
\begin{coro}
Let $C:I\times \mathcal{H}\rightrightarrows \mathcal{H}$ be a
set-valued map with nonempty, closed and convex values satisfying
$(H''_1)$ and $(H''_2)$. Let $f:I\times \mathcal{H}\times
\mathcal{H}\longrightarrow \mathcal{H}$  such that $(H_4)$, $(H_5)$
and $(H_6)$ are satisfied. Then for any $(u_0,v_0)\in
\mathcal{H}\times C(0,u_0)$, there exists a BVRC solution $ (u, v)
:I \longrightarrow \mathcal{H}\times \mathcal{H}$ to  the problem
$$
\begin{cases}
u(0)=u_0, v(0)=v_0;\\
 v(t)\in C(t, u(t))\;\;\;\forall t\in I;\\
\displaystyle\frac{du}{d\nu }(t)=v(t)\;\;\;d\nu-a.e.\,t\in
   I;\\
   -\displaystyle\frac{dv}{d\nu }(t)\in N_{C(t,u(t))} (v(t))+ f(t,u(t), v(t))  \frac{d\lambda} {d\nu}(t)  \;\;\;d\nu-a.e.\,t\in
   I.
\end{cases}
$$
\end{coro}

\end{document}